\documentclass{article}
\usepackage[utf8]{inputenc}
\usepackage[usenames,dvipsnames]{xcolor}
\usepackage{amsmath,amssymb,graphicx,tikz,todonotes}
\usepackage{hyperref}
\usepackage[margin=1in]{geometry}
\usepackage{layouts}
\usepackage[titletoc,title]{appendix} 
\usepackage{tabularx}
\usepackage{placeins}
\usepackage{mathtools}
\usepackage[normalem]{ulem}
\usepackage{xcolor}
\usepackage[style=numeric, citestyle=numeric, backend=biber]{biblatex}

\newcommand{\bi}[1]{\mbox{\boldmath $#1$}}
\def\dV0{{\rm d}V_0}
\def\V0{V_0}

\usepackage[font=small,labelfont=bf]{caption}

\title{Gradient-enhanced crystal plasticity coupled with phase-field fracture modeling}
\author{
Kim Louisa Auth\textsuperscript{a,}\footnote{Corresponding author: \texttt{kim.auth@chalmers.se} (Kim Louisa Auth)}, 
Jim Brouzoulis\textsuperscript{b}, Magnus Ekh\textsuperscript{a}}
\date{}

\begin{document}
\maketitle
\noindent
\textsuperscript{a} \textit{\small Division of Material and Computational Mechanics, Department of Industrial and Materials Science, Chalmers University of Technology, 41296 Sweden}
\newline
\textsuperscript{b} \textit{\small Division of Dynamics, Department of Mechanics and Maritime Sciences, Chalmers University of Technology, 41296 Sweden}

\section*{Abstract}
This study addresses ductile fracture of single grains in metals by modeling of the formation and propagation of transgranular cracks. 
A proposed model integrates gradient extended hardening, phase-field modeling for fracture, and crystal plasticity. It is presented in a thermodynamical framework in large deformation kinematics and accounts for damage irreversibility.
A micromorphic approach for variationally and thermodynamically consistent damage irreversibility is adopted.
The main objective of this work is to analyze the capability of the proposed model to predict transgranular crack propagation. Further, the micromorphic approach for damage irreversibility is evaluated in the context of the presented ductile phase-field model.
This is done by analyzing the impact of gradient-enhanced hardening considering micro-free and micro-hard boundary conditions, studying the effect of the micromorphic regularization parameter, evaluating the performance of the model in ratcheting loading and and testing its capability to predict three-dimensional crack propagation.
In order to solve the fully coupled global and local equation systems, a staggered solution scheme that extends to the local level is presented.
\section{Introduction}
The initiation and propagation of cracks, which play an important role in determining the lifespan of engineering components, are preceded by significant plastic deformation. To simulate and predict the fracture process, particularly the initiation and growth of short (microscopic) cracks evolving into macroscopic cracks, it is imperative to delve into detailed modeling at the grain scale.

Damage and fracture within the grains may occur in the slip planes as a consequence of the localization of plastic slip. The propagation of small cracks follows crystallographic directions \cite{Rovinelli2018}, making it natural to couple damage modeling with crystal plasticity (cf. \cite{Flouriot2003}, \cite{Aslan2011}).
In order to capture size-dependent behavior resulting from geometrical necessary dislocations at grain boundaries (or other obstacles), gradient-extended crystal plasticity models have been introduced, see e.g. \cite{Gurtin2002}, \cite{Evers2004}, \cite{Ekh2007}. These models introduce a length-scale parameter, providing a size-dependent response and helping to overcome mesh-dependence issues for softening behavior \cite{DeBorst1999}.

For modeling the fracture of grains, a combination of damage mechanics and crystal plasticity has been extensively employed, as for example seen in \cite{Aslan2011} and \cite{Ekh2004}. The phase-field approach has in recent years become a popular choice for modeling fracture. In this approach, the discrete crack is represented as a diffuse field whose width is determined by a length scale parameter. The advantageous feature of phase-field modeling in fracture lies in its capability to predict arbitrary crack propagation patterns, including crack branching and coalescence. Moreover, in elasticity, it aligns with classical fracture mechanics, when the length scale parameter approaches zero.
This modeling approach has found application in predicting ductile fracture in several works, as discussed in \cite{Alessi2018} and references therein, and has been coupled with crystal plasticity in \cite{HernandezPadilla2014}, \cite{DeLorenzis2016}, \cite{Maloth2023}. It has been shown to be a promising approach for transgranular fracture of metals.

Phase-field models, however, require special treatment in order to ensure damage irreversibility. While a number of different approaches have been suggested to address its numerical treatment, e.g. \cite{Bourdin2008}, \cite{Miehe2010_thermodynamically_consistent}, \cite{Gerasimov2019}, \cite{Alessi2015}, one of the most commonly adopted approaches is the so called history-variable approach by \citeauthor{Miehe2010_historyvariable} \cite{Miehe2010_historyvariable}.
This approach introduces the history variable in the strong form of the phase field equation. It thereby renders it impossible to retrieve the energy functional back from the strong from and thus results in the loss of variational consistency.
\citeauthor{Bharali2023} \cite{Bharali2023} have recently proposed to use a micromorphic approach instead, offering the advantage that irreversibility can be directly formulated on a local damage field. The micromorphic extension, as presented in \cite{Forest2009}, has in earlier works addressed challenges such as avoiding mesh dependence for softening behavior \cite{Dimitrijevic2011}, model size-dependent hardening \cite{Dimitrijevic2011}, and enhancing numerical robustness \cite{Miehe2017}.

In this contribution, we aim at integrating gradient extended hardening and phase-field modeling within a crystal plasticity framework. The model is presented in a thermodynamical framework with large deformation kinematics. The incorporation of gradient extended hardening is motivated by the resistance to edge dislocation motion at grain boundaries which affects the stress field and therefore the fracture behavior. The incorporation of gradient hardening in a ductile crystal plasticity model is an extension to earlier work in the field. A primary objective of the paper is to analyze model behavior for different choices of gradient hardening boundary conditions and to explore the capability of the micromorphic damage irreversibility approach to interact with a more complex phase-field fracture model than the one presented in \cite{Bharali2023}.

The structure of the paper is as follows: Section~\ref{section:thermodynamics} provides an overview of the thermodynamic framework underlying the model. Section~\ref{section:prototype} introduces a prototype large strain crystal plasticity model, while Section~\ref{section:irreversibility} deals with formulations concerning damage irreversibility. The weak formulation of the balance equations is outlined in Section~\ref{section:weak}. Section~\ref{section:implementation} comments on important details of the numerical implementation, in particular on the staggered solution scheme in presence of multiple local variables. To highlight aspects of the proposed model, numerical results for single crystal boundary value problems in two-dimensional (2D) and three-dimensional (3D) settings are presented and analyzed in Section~\ref{section:examples}. Concluding remarks are given in Section~\ref{section:concluding}.
\section{Thermodynamic modeling framework}
\label{section:thermodynamics}
In this section, a thermodynamic framework for a phase field fracture model based on an underlying gradient crystal plasticity formulation is presented. The derivations follow \citeauthor{Ekh2007} \cite{Ekh2007} for a gradient crystal plasticity model but are here extended with phase-field fracture. For comparison, a derivation based on the microforce balances is given in Appendix~\ref{section:microforce}. 

We formulate the model in a large strain setting and assume a multiplicative decomposition of the deformation gradient $\bi F$ into an elastic part $\bi F_{\rm e}$ and a plastic part $\bi F_{\rm p}$
\begin{equation}
    \bi F=\bi F_{\rm e} \cdot \bi F_{\rm p} \,.
\end{equation}
Further, we introduce isotropic hardening strains $k_\alpha$ on the slip systems $\alpha$ and a damage (phase-field) variable $d$. The free energy  $\Psi$ is then assumed to depend on the elastic Cauchy-Green deformation gradient $\bi C_{\rm e}=\bi F_{\rm e}^{\rm T} \cdot \bi F_{\rm e}$, the state variable $q$,  the isotropic hardening strains $k_\alpha$, the damage variable $d$ and the spatial gradients $\bi \nabla_0 k_\alpha$ and $\bi \nabla_0 d$
\begin{equation}
    \Psi=\Psi\left(\bi C_{\rm e}, \, q, \, k_\alpha, \, \boldsymbol{\nabla}_0 \, k_\alpha, \, d, \, \boldsymbol{\nabla}_0 d \, \right)
    \label{eq:generic_free_energy}
\end{equation}
The dissipation inequality under quasistatic and isothermal conditions is given by
\begin{equation}
    \label{eq:dissipation_inequality_start}
    \int_{\V0} \boldsymbol P: \dot{\bi F} \, \dV0-\int_{\V0} \dot{\Psi} \, \dV0 \geq 0 \,.
\end{equation}
where $\bi P$ is the first Piola-Kirchhoff stress, $\bi F$ is the deformation gradient and $\V0$ represents the initial domain with boundary $\Gamma_0$.
Introducing the free energy (\ref{eq:generic_free_energy}) into the dissipation inequality (\ref{eq:dissipation_inequality_start}) and using the standard Coleman-Noll arguments \cite{Coleman1963} 
yields the elastic second Piola-Kirchhoff stress $\bi S_{\rm e}$ as (see e.g. \cite{Simo1988})
\begin{equation}
    \bi S_{\rm e}=2 \, \frac{\partial \Psi}{\partial \bi C_{\rm e}} \,,
    \label{eq:elastic_2nd_PK_stress}
\end{equation}
and we obtain the reduced dissipation inequality
\begin{equation}
\label{eq:micromorphic_gradientplastic_dissipation_inequality}
\mathcal{D} =
        \int_{V_0}
        \left(
      {\boldsymbol{M}}_\mathrm{e}  
                :
                {\boldsymbol{L}}_\mathrm{p}  
                -
        \frac{\partial\Psi}{\partial q} \, \dot{q}
    -
        \frac{\partial\Psi}{\partial k_\alpha} \, \dot{k}_\alpha
    -
        \frac{\partial\Psi}{\partial \boldsymbol{\nabla}_0 k_\alpha} \cdot \boldsymbol{\nabla}_0 \dot{k}_\alpha
    -
        \frac{\partial\Psi}{\partial d} \, \dot{d}
    -
        \frac{\partial\Psi}{\partial \boldsymbol{\nabla}_0 d} \cdot \boldsymbol{\nabla}_0 \dot{d}
    \right)
    \; \mathrm{d}V_0
    \geq
    0
    \,.
\end{equation}
In equations (\ref{eq:elastic_2nd_PK_stress})-(\ref{eq:micromorphic_gradientplastic_dissipation_inequality}), 
the Mandel stress ${\boldsymbol{M}}_\mathrm{e}=\bi C_{\rm e} \cdot \bi S_{\rm e}$, like the second Piola-Kirchhoff stress on the intermediate configuration, as well as the plastic velocity gradient ${\boldsymbol{L}}_\mathrm{p}=\dot{\bi F}_{\rm p} \cdot {\bi F}_{\rm p}^{-1}$ were introduced.
Using the divergence theorem, the reduced dissipation inequality can be rewritten as
\begin{equation}
    \mathcal{D} =
        \int_{V_0}
                {\boldsymbol{M}}_\mathrm{e}  
                :
                {\boldsymbol{L}}_\mathrm{p}
                + Q \, \dot{q}
   + \kappa_\alpha
            \, \dot{k}_\alpha +
           Y\, \dot{d}
        \; \mathrm{d}V_0 
    +
        \int_{\Gamma_0}
            \kappa_\alpha^{\Gamma}  \, \dot{k}_\alpha
       +Y^{\Gamma}
            \, \dot{d}
        \, \mathrm{d}V_0
    \geq
    0
    \,,
    \label{eq:diss_inequality_2}
\end{equation}
where $Q=-\partial \Psi/\partial q$ and the gradient extended dissipative hardening stress $\kappa_\alpha$ and its boundary "traction" $\kappa_\alpha^\Gamma$ were introduced as
\begin{align}
      \kappa_\alpha&=-\frac{\partial\Psi}{\partial k_\alpha} +\boldsymbol{\nabla}_0 \cdot \frac{\partial\Psi}{\partial\boldsymbol{\nabla}_0 k_\alpha} \,,
      \label{eq:kappa_def}
      \\
      \kappa_\alpha^{\Gamma}&= \boldsymbol{N}
            \cdot
            \frac{\partial\Psi}{\partial\boldsymbol{\nabla}_0 k_\alpha}
        \,,
\end{align}
where $\bi N$ is the unit normal to $\Gamma_0$. The corresponding dissipative quantities for the phase field are
\begin{eqnarray}
      Y &=& -\frac{\partial\Psi}{\partial d} 
                + 
                \boldsymbol{\nabla}_0 \cdot \frac{\partial\Psi}{\partial\boldsymbol{\nabla}_0 d}
               \label{eq:Y_def}
                \\
       Y^{\Gamma} &=&   \boldsymbol{N}
            \cdot
            \frac{\partial\Psi}{\partial\boldsymbol{\nabla}_0 d}       
\end{eqnarray}
The phase field equation is obtained by assuming that $Y=0$ (purely energetic) 
\begin{equation}
   -\frac{\partial\Psi}{\partial d} 
                + 
                \boldsymbol{\nabla}_0 \cdot \frac{\partial\Psi}{\partial\boldsymbol{\nabla}_0 d}=0 
\label{eq:pf_base}
\end{equation}
and thereby it can be noted that the phase-field will not contribute to the dissipation on $V_0$.
Hence, the reduced dissipation inequality becomes
\begin{equation}
    \mathcal{D} =
        \int_{V_0}
        \left(
                {\boldsymbol{M}}_\mathrm{e}  
                :
                {\boldsymbol{L}}_\mathrm{p}
   +  Q \, \dot{q} +\kappa_\alpha
            \, \dot{k}_\alpha 
        \right)
        \; \mathrm{d}V_0 
    +
        \int_{\Gamma_0}
        \left(
            \kappa_\alpha^{\Gamma}  \, \dot{k}_\alpha
       +Y^{\Gamma}
            \, \dot{d}
        \right)
        \, \mathrm{d}V_0
    \geq
    0
\label{eq:reduced_dissipation_inequality}
\end{equation}
By comparing with the derivation in Appendix~\ref{section:microforce}, we can conclude the dissipation inequality is not defined locally and that also boundary terms are obtained. 
These must be checked when formulating boundary conditions for the field equations. The resulting field equations (\ref{eq:kappa_def}) and (\ref{eq:pf_base}) are the same.  

\section{Prototype crystal plasticity model}
\label{section:prototype}
We will assume the following form of the free energy
\begin{equation}
    \Psi
    =
    g_\mathrm{e} \left(d, \, \epsilon^\mathrm{p}\right) \, \hat{\Psi}_\mathrm{e}\left(\boldsymbol{C}_\mathrm{e} \right)
    +
    g_\mathrm{p}\left(d\right) 
    \sum_\alpha{
        \hat{\Psi}_\mathrm{p}\left(k_\alpha, \boldsymbol{\nabla}_0 k_\alpha\right)
    }
    +
    {\Psi}_\mathrm{d}\left(d, \boldsymbol{\nabla}_0 d\right)
\end{equation}
The effective (undamaged) elastic part of the free energy $\hat{\Psi}_\mathrm{e}$ is assumed to be of Neo-Hookean type
\begin{equation}
  \hat{\Psi}_\mathrm{e}\left(\boldsymbol{C}_\mathrm{e} \right)=
  \frac{\mu}{2} \, \left( {\rm tr}(\boldsymbol{C}_\mathrm{e}) -3 \right)- \mu \, {\rm ln}(J_{\rm e})+\frac{\lambda}{2} \, \left( {\rm ln}(J_{\rm e}) \right)^2
  \quad \mbox{with} \; J_{\rm e}^2={\rm det}(\bi C_{\rm e})
\end{equation}
where $\lambda$ and $\mu$ are the elastic Lamé constants. For simplicity, we have assumed elastic isotropy and disregarded a tension-compression split.  
The elastic degradation function $g_\mathrm{e}$ is adopted from the ductile fracture model presented by \citeauthor{Ambati2015} \cite{Ambati2015}
\begin{equation}
\label{eq:degradation_function}
    g_\mathrm{e}\left(d, \, \epsilon^\mathrm{p}\right)
    = 
    \left(1 - d\right)^{2\left(\epsilon^\mathrm{p} / \epsilon^\mathrm{p}_\mathrm{crit}\right)^n}
    \,,
\end{equation}
where $\epsilon^\mathrm{p}$ is the accumulated plastic strain (in the thermodynamic modeling framework above represented by $q$). It can be noted that the degradation is only active when $\epsilon_{\rm p}>0$ and that the parameters $\epsilon^\mathrm{p}_\mathrm{crit}$ and $n>0$ control how the degradation increases when $\epsilon^\mathrm{p}$ increases.  As shown by \cite{Ambati2015}, the formulation gives a positive contribution to the dissipation, i.e. $Q \ \dot{q} \geq 0$. 

The effective plastic free energy $\hat{\Psi}_\mathrm{p}$ is is chosen as (compare \cite{Ekh2007})
\begin{equation}
\label{eq:plastic_free_energy}
 \hat{\Psi}_{\rm p}=\frac{1}{2} \, \sum_\alpha H_\alpha \, k_\alpha^2+
\frac{l_g^2}{2} \, \sum_\alpha H^g_\alpha \, (\bar{\bi s}_\alpha \cdot \nabla_0 k_\alpha)^2  \,, 
\end{equation}
where $H_\alpha$ is the isotropic hardening modulus, $H^\mathrm{g}_\alpha$ is the gradient-enhanced hardening modulus and $l^\mathrm{g}$ is the length scale for gradient-enhanced hardening. 
We adopt the standard assumption that the slip direction $\bar{\bi s}_\alpha$ and normal vector to the slip plane $\bar{\bi m}_\alpha $ on the intermediate  configuration are fixed (and equal to their corresponding vectors on the undeformed configuration).
For simplicity, the plastic degradation function is chosen as $g_{\rm p}=1$.
The yield function $\Phi_\alpha$ is defined in terms of the effective Schmid stress $\hat{\tau}_\alpha$ as
\begin{eqnarray}
   \Phi_\alpha
   =
   \left| \hat{\tau}_\alpha \right|
   - ( \tau_{\rm y}+\kappa_\alpha )
\end{eqnarray}
with $\hat{\tau}_\alpha={\tau}_\alpha\,/\,g_\mathrm{e}\left(d, \, \epsilon^\mathrm{p}\right)$ wherein $\tau_\alpha$ is the standard crystal plasticity Schmid stress $\tau_\alpha=\bi{M}_{\rm e}: \left( \bar{\bi s}_\alpha \otimes \bar{\bi m}_\alpha \right)$.  Furthermore, $\tau_{\rm y}$ is the initial yield stress.
The evolution equation for the plastic velocity gradient is assumed to be of associative type 
\begin{equation}
\label{eq:evolution_Fp}
    \bar{\boldsymbol{L}}_{\rm p}=
    \dot{\bi F}_{\rm p} \cdot \bi F_{\rm p}^{-1}= \sum_{\alpha=1}^{N_\alpha} \dot{\lambda}_\alpha \, \frac{\partial \Phi_\alpha}{\partial {{\bi M}}_{\rm e}}
    = \sum_{\alpha=1}^{N_\alpha} \frac{\dot{\lambda}_\alpha}{g_\mathrm{e}\left(d, \, \epsilon^\mathrm{p}\right)} \, \left(
\bar{\boldsymbol{s}}_\alpha \otimes \bar{\boldsymbol{m}}_\alpha
\right)
\end{equation}
and we apply a viscoplastic regularization for the multiplier
\begin{equation}
    \dot{\lambda}_\alpha = \frac{1}{ t^\ast}\, \left< \frac{\Phi_\alpha}{\sigma_\text{d}} \right> ^ {m} \,,
\end{equation}
where $t^\ast$, $m$ and $\sigma_\text{d}$ control the viscosity of the model and $<\bullet>$ denotes Macualey brackets. The accumulated plastic strain $\epsilon^{\rm p}$ is based on $\dot{\lambda}_\alpha$ defined as
\begin{equation}
\label{eq:evolution_ep}
    \epsilon_{\rm p}
    =
    \int_0^t 
    \sqrt{
        \sum_{\alpha=1}^{N_\alpha} \dot{\lambda}_\alpha^2
    }
    \, {\rm d}t \,.
\end{equation}
The gradient extended hardening stress $\kappa_\alpha$ is derived from Equation (\ref{eq:kappa_alpha}) as
\begin{equation}
 \kappa_\alpha=  -H_\alpha \, k_\alpha+H^g_\alpha \, l_g^2 \, \bar{\bi s}_\alpha \cdot (
\boldsymbol{\nabla}_0 \otimes \boldsymbol{\nabla}_0  k_\alpha)
\cdot     \bar{\bi s}_\alpha  
\label{eq:kappa_alpha}
\end{equation}
and the evolution of the hardening strains $k_\alpha$ is also assumed to be of associative type
\begin{eqnarray}
\label{eq:evolution_k}
    \dot{k}_\alpha=\dot{\lambda}_\alpha \, \frac{\partial \Phi_\alpha}{\partial \kappa_\alpha}=-\dot{\lambda}_\alpha \,.
\end{eqnarray}
These assumptions for the hardening can be extended, see e.g. \cite{Bargmann2010} to account for more complex models such as kinematic hardening, cross-hardening and nonlinear hardening.

For the phase-field fracture model, the free energy contribution is based on an AT2 surface energy functional $\Gamma_\mathrm{d}$, cf. \cite{Ambrosio1990} 
\begin{equation}
 {\Psi}_\mathrm{d}\left(d, \boldsymbol{\nabla}_0 d\right)=
 \mathcal{G}^\mathrm{d}_0 \, \Gamma_\mathrm{d}\left(d, \boldsymbol{\nabla}_0 d \right) \,
 \quad
 \mbox{with} \;
   \Gamma_\mathrm{d}  =
    \frac{1}{2 \, \ell_0}
    \left(
        d^2 + \ell_0^2 \left| \boldsymbol{\nabla}_0 d \right|^2
    \right)
\end{equation}
where $\mathcal{G}^\mathrm{d}_0$ represents fracture toughness and $\ell_0$ is the length-scale parameter controlling the width of the diffuse crack model.
By inserting the choices of degradation functions and ${\Psi}_\mathrm{d}$,
the phase field equation (\ref{eq:pf_base}) becomes
\begin{eqnarray}
2 \, \left(\epsilon^\mathrm{p} / \epsilon^\mathrm{p}_\mathrm{crit}\right)^n \,(1-d)^{2 \, \left(\epsilon^\mathrm{p} / \epsilon^\mathrm{p}_\mathrm{crit}\right)^n-1}
\hat{\Psi}_\mathrm{e}
        - \frac{\mathcal{G}_0^\mathrm{d}}{\ell_0} \,d 
    + \mathcal{G}_0^\mathrm{d}\,\ell_0 \,\boldsymbol{\nabla}_0 \cdot \boldsymbol{\nabla}_0 d
    = 0 
\label{eq:pf_prototype}
\end{eqnarray}
which is similar to the formulation used in \cite{DeLorenzis2016} but extended with the exponent $n$. 

\section{Irreversibility}
\label{section:irreversibility}
As discussed in the introduction, one of the most common approaches for enforcing damage irreversibilty is the  history variable approach introduced by \citeauthor{Miehe2010_historyvariable} \cite{Miehe2010_historyvariable}
where the effective elastic free energy $\hat{\Psi}_\mathrm{e}$ in the phase-field equation (Equation (\ref{eq:pf_prototype})) 
is replaced by a history variable
\begin{equation}
 \mathcal{H}(t)=\max_{\tilde{t} \leq t} \hat{\Psi}_{\rm e} \,.
\end{equation}
This approach though has been shown to be variationally inconsistent and it is not trivial to show its thermodynamic consistency upon unloading. Within this work we therefore explore a micromorphic approach \cite{Forest2009}, which has recently been shown to allow for a variationally consistent framework for locally enforced damage irreversibilty \cite{Bharali2023}.
It introduces an additional local variable $\varphi$ and includes a penalty term in the free energy that connects the global damage $d$ with the (new) local damage $\varphi$. The local damage $\varphi$ then replaces $d$ in all terms except for the new penalty term and the damage gradient term.
\begin{equation}
    \label{eq:free_energy_micromorphic}
    \Psi
    =
    g_\mathrm{e} \left(\varphi, \, \epsilon^\mathrm{p}\right) \, \hat{\Psi}_\mathrm{e}\left(\boldsymbol{C}_\mathrm{e} \right)
    +
    g_\mathrm{p}\left(\varphi\right) 
    \sum_\alpha{
        \hat{\Psi}_\mathrm{p}\left(k_\alpha, \boldsymbol{\nabla}_0 k_\alpha\right)
    }
    +
    {\Psi}_\mathrm{d}\left(\varphi, \boldsymbol{\nabla}_0 d\right)
    +
    \frac{\alpha}{2} \left( \varphi-d\right)^2
\end{equation}
The dissipation inequality (\ref{eq:diss_inequality_2}) is thereby modified to
\begin{equation}
    \mathcal{D} =
        \int_{V_0}
                {\boldsymbol{M}}_\mathrm{e}  
                :
                {\boldsymbol{L}}_\mathrm{p}
                + Q \, \dot{q}
   + \kappa_\alpha
            \, \dot{k}_\alpha +Y_\varphi\, \dot{\varphi}+
           Y_\mathrm{d}\, \dot{d}
        \; \mathrm{d}V_0 
    +
        \int_{\Gamma_0}
            \kappa_\alpha^{\Gamma}  \, \dot{k}_\alpha
       +Y^{\Gamma}
            \, \dot{d}
        \, \mathrm{d}V_0
    \geq
    0
    \label{eq:diss_inequality_2b}
\end{equation}
where equivalently to the procedure leading to Equation (\ref{eq:pf_base}), $Y_\mathrm{d}=0$ yields the global phase-field equation
\begin{equation}
    Y_\mathrm{d}
    =
    \alpha \, (\varphi-d)
    + \mathcal{G}_0^\mathrm{d}\,\ell_0 \,\boldsymbol{\nabla}_0 \cdot \boldsymbol{\nabla}_0 d
    = 0 \,
\label{eq:pf-micromorphic}
\end{equation}
and hence the global phase field $d$ does not contribute to the dissipation on $V_0$.
The micromorphic approach introduces an additional regularization to the model. In \cite{Miehe2017} it was mainly used for robustness of numerical implementation but can as suggested by \cite{Bharali2023} conveniently be used to ensure irreversibility.
The evolution of the local phase field $\varphi$ is then derived from the inequality $Y_\varphi \, \dot{\varphi} \geq 0$. For the suggested choice of model $Y_\varphi$ is given by
\begin{equation}
\label{eq:local_residual_pf}
Y_\varphi=- \frac{\partial \Psi}{\partial \varphi}=
2 \, \left(\epsilon^\mathrm{p} / \epsilon^\mathrm{p}_\mathrm{crit}\right)^n \,(1-\varphi)^{2 \, \left(\epsilon^\mathrm{p} / \epsilon^\mathrm{p}_\mathrm{crit}\right)^n-1}
\hat{\Psi}_\mathrm{e}
        - \frac{\mathcal{G}_0^\mathrm{d}}{\ell_0} \,\varphi 
-\alpha \, (\varphi-d) \,.
\end{equation}
For pure loading, the local phase field $\varphi$ can be computed from assuming that $Y_\varphi$ is energetic, i.e. $Y_\varphi=0$. In order to obtain a thermodynamically consistent formulation for unloading however, the full inequality must be considered. By introducing an intermediate local variable $\tilde{\varphi}=\arg\left\{Y_\varphi(\varphi)=0\right\}$, the  Karush-Kuhn-
Tucker conditions to ensure irreversibility of $\varphi$ can be formulated as
\begin{equation}
    \label{eq:varphi_KKT}
    \dot{\varphi} \geq 0 \,, \quad 
	    \dot{\varphi} \, f_\mathrm{\varphi}=0 \,, \quad  
	    f_\mathrm{\varphi} = \tilde{\varphi}-\varphi \leq 0 \,.
\end{equation}
Thereby, we obtain $Y_\varphi \, \dot{\varphi} =0$, since $Y_\varphi=0$ during loading ($\varphi=\tilde{\varphi}$) and $\dot{\varphi}=0$ during unloading ($\varphi=\mathrm{const}$). Notice that the local phase field here becomes a history variable.

\section{Weak form of balance equations}
\label{section:weak}
The weak form of balance of momentum when neglecting inertial forces and body forces is expressed in terms of the first Piola-Kirchhoff stress $\boldsymbol{P}$ as
\begin{equation}
    \mathcal{R}^\mathrm{u}
    =
   \int_{V_0}
        \boldsymbol{P} : \left( \delta \boldsymbol{u} \otimes 
        \boldsymbol{\nabla}_0 \right)
        \,
    \mathrm{d}V_0
    -\int_{\partial V_0}
        \boldsymbol{t}_0^\ast
        \cdot \delta \boldsymbol{u}
        \,
    \mathrm{d}A_0
    = 0 \,,
\label{eq:weak_equil}
\end{equation}
where  $\boldsymbol{t}_0^\ast$ is a prescribed traction on the boundary $\partial V_0$. 
The expression for the gradient extended hardening stress $\kappa_\alpha$ in Equation (\ref{eq:kappa_alpha}) is, due to the spatial gradients, also a field equation. We adopt the dual mixed procedure described in \cite{Svedberg1998} and introduce
\begin{equation}
    \boldsymbol{g}_\alpha = \boldsymbol{\nabla}_0 k_\alpha 
    \label{eq:strong_galpha}
\end{equation}
whereby
\begin{eqnarray}
    \kappa_\alpha
    =
    -H_\alpha \, k_\alpha
    +
    H^g_\alpha \, l_g^2 \, \bar{\bi s}_\alpha
    \cdot 
    \boldsymbol{\nabla}_0 \boldsymbol{g}_\alpha   
    \cdot
    \bar{\bi s}_\alpha
\end{eqnarray}
becomes a local equation. 
Instead, the weak form of Equation (\ref{eq:strong_galpha}) is introduced as a field equation and by using the divergence theorem we obtain
\begin{equation}
    \mathcal{R}^g
    =
    \int_{\V0} 
         \boldsymbol{g}_\alpha \, \delta  \boldsymbol{g}_\alpha
    \, \dV0
     -
     \int_{\partial V_0} k_\alpha \, \bi N \cdot \bar{\bi s}_\alpha \, \delta g_\alpha \, {\rm d}A_0
     +
     \int_{\V0} k_\alpha \, \boldsymbol{\nabla}_0 ( \delta  \bi g_\alpha) \cdot  \bar{\bi s}_\alpha \, \dV0    
     =
     0 \,.
\label{eq:weak_galpha}
\end{equation}
This procedure has been shown to be numerically robust (\cite{Ekh2007}, \cite{Carlsson2017}) and has the same benefit as a micromorphic approach that no special treatment to account for loading/unloading conditions is needed since it can be controlled locally in the material points.
Finally, the phase-field equation (\ref{eq:pf-micromorphic}) is formulated in weak form as
\begin{equation}
\mathcal{R}^\mathrm{d}
=
\int_{V_0}
        \alpha \left(\varphi - d\right) \, \delta d \,
    \mathrm{d}V_0
    +
     \int_{\partial V_0}
        {\,\mathcal{G}_0^\mathrm{d}\,\ell_0} \, \bi{N} \cdot 
        \boldsymbol{\nabla}_0 d \,  \delta d \,
    \mathrm{d}A_0
    -
    \int_{V_0}
        {\,\mathcal{G}_0^\mathrm{d}\,\ell_0} \, k        \boldsymbol{\nabla}_0 d \cdot \boldsymbol{\nabla}_0 \delta d \,
    \mathrm{d}V_0
    = 0
\label{eq:weak_pf}
\end{equation}
where the standard boundary condition $\bi{N} \cdot 
\boldsymbol{\nabla}_0 d =0$ will be assumed whereby the boundary integral term disappears.

\section{Numerical implementation}
\label{section:implementation}
\begin{figure}[htb]
    \centering
    \includegraphics{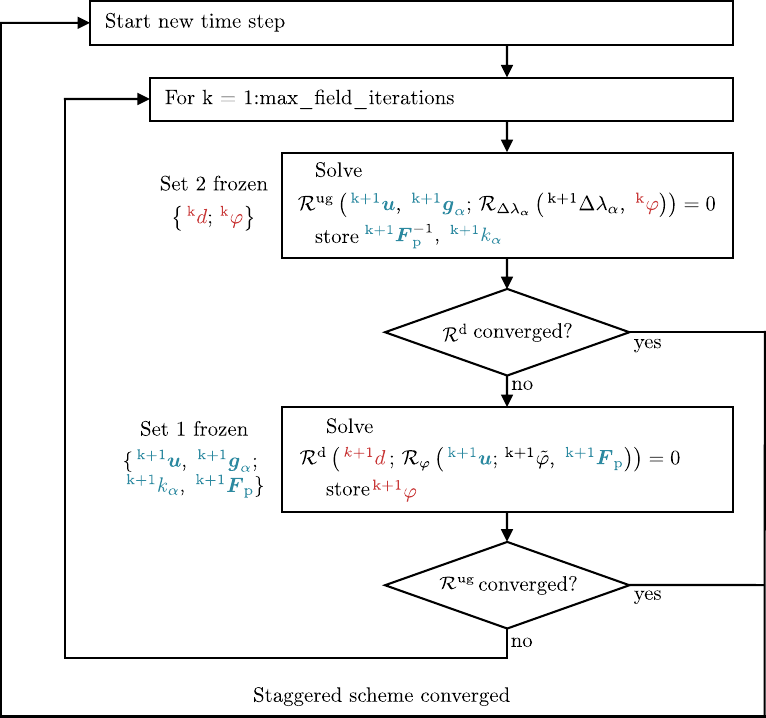}
    \caption{Staggered iteration scheme for the global system solver. The global and local variables are grouped in two groups extending the staggered solution scheme to the local equation system. The staggered algorithm solves for one set of variables at a time while freezing the other set. The iteration between the sets is stopped once the global residuals fulfill the all convergence criteria presented in Table \ref{tab:global_tolerances} at the same time.}
    \label{fig:staggered_scheme}
\end{figure}

Time integration schemes are needed for the plastic evolution equations, Equation (\ref{eq:evolution_Fp}) and (\ref{eq:evolution_k}), as well as for the accumulated plastic strain, Equation (\ref{eq:evolution_ep}) and the evolution of the local phase field, Equation (\ref{eq:varphi_KKT}).
Backward Euler time integration is applied to the evolution equations for the plastic deformation gradient $\boldsymbol{F}_\mathrm{p}$ and the hardening strains $k_\alpha$.
For time points $\,^\mathrm{n+1}t$ and $\,^\mathrm{n}t$, a time step $\Delta t = \,^\mathrm{n+1}t - \,^\mathrm{n}t$ and $\Delta \lambda_\alpha = \Delta t \,^\mathrm{n+1}\dot{\lambda}_\alpha$, we obtain the following expressions
%For a new time step $n+1$, a previous time step $n$ and $\Delta \lambda_\alpha = (\,^\mathrm{n+1}t - \,^\mathrm{n}t)\,^\mathrm{n+1}\dot{\lambda}_\alpha$ we obtain the following expressions
\begin{align}
    \,^\mathrm{n+1}k_\alpha
    &=  
    \,^\mathrm{n} k_\alpha - \,^\mathrm{n+1} \Delta\lambda_\alpha 
    \\
    \,^\mathrm{n+1}\boldsymbol{F}_\mathrm{p}^{-1} 
    &=
    \,^\mathrm{n} \boldsymbol{F}_\mathrm{p}^{-1}
    \cdot
    \left(
        \boldsymbol{I} - 
        \sum_\alpha
            \frac{\,^\mathrm{n+1}\Delta\lambda_\alpha}{g_\mathrm{e}(\,^\mathrm{n+1}\varphi)}
            \left(\bar{\boldsymbol{s}}_\alpha \otimes \bar{\boldsymbol{m}}_\alpha \right)
            \,
            \mathrm{sign}\left(\,^\mathrm{n+1}\tau_\alpha\right)
    \right) \,.
\end{align}
The signs of the Schmid stresses $\,^\mathrm{n+1}\tau_\alpha$ are computed based on the elastic trial stress.
The evolution of the local phase field $\varphi$ is discretized as
\begin{equation}
    \,^\mathrm{n+1}\varphi
    =
    \max\left(\,^\mathrm{n+1}\tilde{\varphi},\, \,^\mathrm{n}\varphi\right)
\end{equation}
and thereby accounting for the irreversibility condition. An explicit scheme is applied to the accumulated plastic strain $\epsilon^\mathrm{p}$, such that the elastic degradation is computed based on the accumulated plastic strain from the previous time step $\,^\mathrm{n}\epsilon^\mathrm{p}$
\begin{equation}
    g_\mathrm{e}\left(\varphi \right)
    = 
    \left(1 - \varphi\right)^{2\left(\,^\mathrm{n} \epsilon^\mathrm{p} / \epsilon^\mathrm{p}_\mathrm{crit}\right)^n}
    \,.
\end{equation}
Thereby, two coupled local residual equations with unknowns $\,^\mathrm{n+1}\tilde{\varphi}$ and $\,^\mathrm{n+1}\Delta\lambda_\alpha$ need to be solved for every time step within the global residual equations (\ref{eq:weak_equil}), (\ref{eq:weak_galpha}) and (\ref{eq:weak_pf}).
\begin{align}
    \mathcal{R}_{\Delta\lambda_\alpha}
    \left(
        \Delta\lambda_\alpha,\, \varphi
    \right)
    &=
    \Delta\lambda_\alpha
    -
    \frac{\Delta t}{t^\ast} \, \left< \frac{\Phi_\alpha
        \left(
            \boldsymbol{C},\, \boldsymbol{\nabla}_0\boldsymbol{g}_\alpha,\,\Delta\lambda_\alpha,\, \varphi
        \right)
    }{\sigma_\mathrm{d}} \right>^\mathrm{m}
    \\
    \mathcal{R}_\varphi
    \left(
        \Delta\lambda_\alpha,\, \tilde{\varphi}
    \right)
        &=
        - \frac{\partial g_\mathrm{e}\left(\tilde{\varphi}\right)}{\partial \tilde{\varphi}}
        \hat{\Psi}_\mathrm{e}
        \left(
            \boldsymbol{C},\, \Delta\lambda_\alpha,\, \tilde{\varphi}
        \right)
        - \frac{\mathcal{G}_0^\mathrm{d}}{\ell_0} \,\tilde{\varphi}
-\alpha \, (\tilde{\varphi}-d)
\end{align}
We apply a staggered algorithm in order to solve the coupled problem. The algorithm is schematically represented in Figure \ref{fig:staggered_scheme}. The equilibrium equation $\mathcal{R}^\mathrm{u}$ and the gradient extended hardening field equation $\mathcal{R}^\mathrm{g}$ are solved in a monolithic way (refered to as $\mathcal{R}^\mathrm{ug}$ in Figure \ref{fig:staggered_scheme}) in the first staggered step and the global phase-field equation $\mathcal{R}^\mathrm{d}$ is solved in a second staggered step. Extending the staggered approach to the local equations, we group variables in two sets corresponding to the staggered steps:
Set 1: \{$\boldsymbol{u}$, $\boldsymbol{g}_\alpha$, $k_\alpha$, $\boldsymbol{F}_\mathrm{p}$\} and set 2: \{$d$, $\varphi$\}. Within each staggered iteration only the corresponding set of variables is updated, while the other set is frozen. Thus, a decoupling of the local equations is achieved and only one of the local residual equations needs to be solved within each staggered step.
This is a similar staggered approach as the one taken by \citeauthor{Ambati2015} \cite{Ambati2015}.

The local equations are solved by Newton iterations up to a tolerance of $10^{-8}$. Field-wise convergence criteria on the global fields and their respective residuals are applied on the global system. A field is considered converged when either the norm of the residuals is below the respective residual tolerance or the norm of the update of field values in a Newton iteration is below the respective field tolerance. 
The tolerances applied for the numerical examples are displayed in Table \ref{tab:global_tolerances}.
\begin{table}
    \centering
    \begin{tabular}{||l c l l||}
        \hline
         ~ & \vtop{\hbox{\strut Global residual}\hbox{\strut equation}}& \vtop{\hbox{\strut Residual}\hbox{\strut tolerance}}& \vtop{\hbox{\strut Field }\hbox{\strut tolerance}} \\ [0.5ex]
        \hline \hline
        Displacement field $\boldsymbol{u}$ & $\mathcal{R}^\mathrm{u}$, Equation \ref{eq:weak_equil} & $10^{-3}\,\mathrm{Nmm}$ & $10^{-6}\,\mathrm{mm}$\\
        \hline
        Hardening strain gradient $\boldsymbol{g}_\alpha$ & $\mathcal{R}^\mathrm{g}$, Equation \ref{eq:weak_galpha} & $10^{-7}\,\mathrm{mm}$& $10^{-8}\,\mathrm{mm}^{-1}$\\
        \hline
        Global phase field $d$ & $\mathcal{R}^\mathrm{d}$, Equation \ref{eq:weak_pf} & $10^{-4}\,\mathrm{Nmm}$& $10^{-8}\,-$\\
        \hline
    \end{tabular}
    \caption{Field-wise tolerances applied for global convergence examination. A field is considered converged when either of the residual or the field tolerance (relating to the update of field variables during a single Newton-correction) is fulfilled.}
    \label{tab:global_tolerances}
\end{table}

\section{Numerical experiments on single crystals}
\label{section:examples}
In this section, the behavior of the proposed prototype model is investigated.
We first demonstrate that the presented model can reproduce well known behavior of crystal plasticity phase field models from the literature. Then we display the effect of added gradient-enhanced hardening.
Subsequently, it is demonstrated that the model is capable of capturing irreversible unloading. Finally, an example of damage development in a three-dimensional setting with inhomogeneities is shown.

\begin{figure}[htb]
    \centering
    \includegraphics{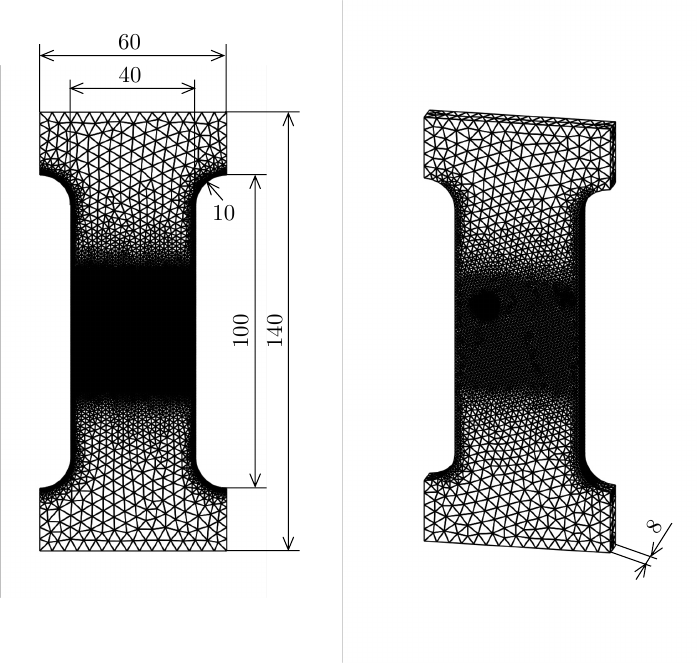}
    \caption{Meshes for the 2D (left) and 3D (right) numerical examples. The cross-sectional geometry of the 3D-example is the same as for the 2D-example and geometrical measures are given in mm. Both examples employ unstructured meshes with a background element size of $4\,\mathrm{mm}$ and mesh refinements at the sides and in the center of the web, where the sample is expected to break. The refinements in the 2D-mesh consist of elements with an average size of $0.5\,\mathrm{mm}$, resulting in $19\,390$ triangular elements in total. In the 3D-mesh, center and sides employ elements of $1\,\mathrm{mm}$ average size and two positions are additionally refined with elements of $0.5\,\mathrm{mm}$ on average for imposing initial material inhomogenities. The 3D mesh consists of $123\,183$ tetrahedral elements in total.}
    \label{fig:meshes}
\end{figure}

The numerical examples, inspired by \cite{DeLorenzis2016}, employ an I-shaped specimen in 2D (plane strain) and 3D. A base set of material parameters is shown in Table \ref{tab:material_parameters}, deviating material parameters are given in the descriptions of the respective numerical examples. Figure \ref{fig:meshes} shows the geometry and meshes for the numerical examples. Both cases represent the same cross-sectional geometry and employ unstructured meshes with a background element size of $4\,\mathrm{mm}$.
Mesh refinements are conducted at the sides and in the center of the web, where the specimen is expected to break. The refinements in the 2D-mesh consist of elements with an average size of $0.5\,\mathrm{mm}$, resulting in $19\,390$ triangular elements in total. In the 3D-mesh, center and sides are meshed with elements of $1\,\mathrm{mm}$ average size and two positions are additionally refined with elements of $0.5\,\mathrm{mm}$ on average for imposing initial material inhomogeneities. The 3D mesh consists of $123\,183$ tetrahedral elements in total.
In all cases, FCC slip systems whose unit cell is aligned with the coordinate axes are used. The slip systems are shown in Table \ref{tab:slipsystems} in Appendix \ref{section:appendix_parameters}. The slip systems with slip directions perpendicular to the loading direction are omitted in all examples, since they have a negligible impact on the plastification. Sections \ref{subsection:example_gradient_hardening} and \ref{subsection:example_micromorphic} thereby use 8 slip systems each. Sections \ref{subsection:example_cyclic} and \ref{subsection:example_3d} only employ the four most active slip systems, that is slip systems $1,\,4,\,7\,\text{and}\,10$.

The I-shaped specimen is clamped on the bottom side and loaded by prescribing a displacement of $10\,\mathrm{mm}/\mathrm{s}$ in the vertical direction on the upper side, while the horizontal displacement is zero.

\subsection{Boundary conditions for gradient-enhanced hardening}
\label{subsection:example_gradient_hardening}
\begin{figure}
    \centering
    \includegraphics{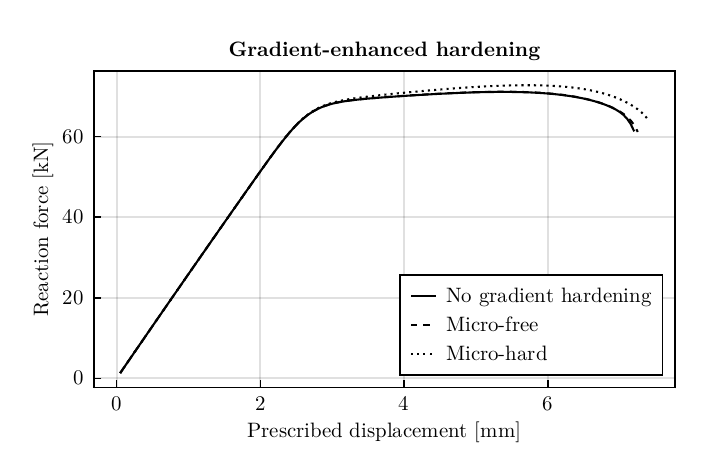}
    \caption{Reaction force response for the I-shaped beam exposed to different boundary conditions, as well as without gradient hardening ($l^\mathrm{g}=0.0\,\text{mm}$). All scenarios reach the softening regime. The micro-hard restriction of slip transfer on the boundaries leads to a stiffer response in the hardening regime. The last point of the curves corresponds to the last converged time step of the respective simulations.}
    \label{fig:reactionforces_gradient_hardening}
\end{figure}
\begin{figure}[htb]
    \centering\includegraphics{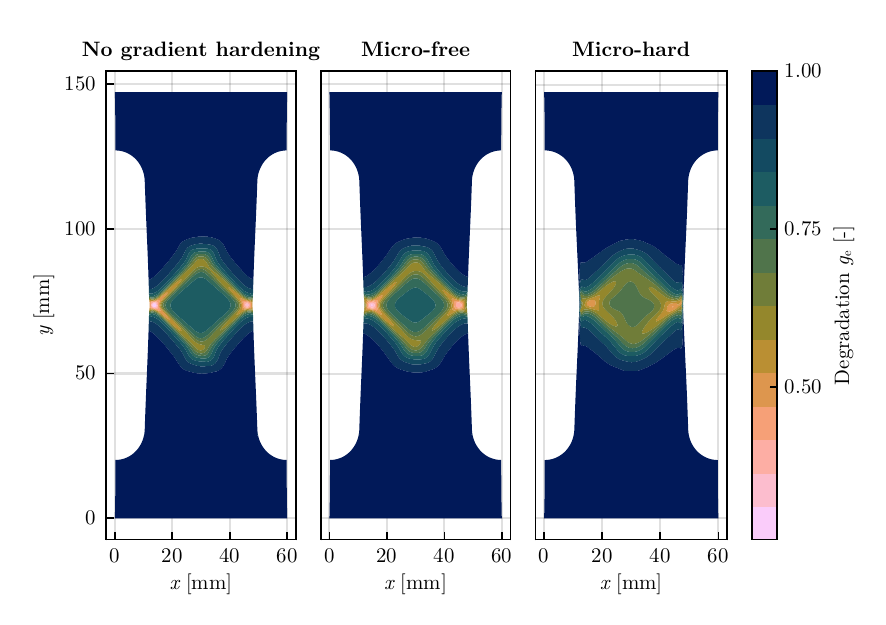}
    \caption{Degradation $g_\mathrm{e}$ at the final step of the three base scenarios. The well known diamond shape crack pattern (compare e.g. \cite{DeLorenzis2016}) is recovered without gradient hardening. The behavior under micro-free boundary conditions is similar, but shows some smoothing of the crack shape. Under micro-hard boundary conditions plastic strains cannot develop on the boundary, thereby preventing the development of damage on the boundary.}
    \label{fig:degradation_gradient_hardening}
\end{figure}
In a first step, we consider the I-shaped specimen under plane strain conditions. We investigate the model behavior without gradient hardening, $l^\mathrm{g}=0.0\,\text{mm}$, as well as for the two trivial possibilities of boundary conditions on the strain gradient fields: micro-hard, $k_\alpha = 0$, and micro-free, $\boldsymbol{g}_\alpha \cdot \bar{\boldsymbol{s}}_\alpha = 0$, slip transfer. Note that both boundary conditions give no dissipation on the boundary, $\kappa_\alpha^{\Gamma}  \, \dot{k}_\alpha=0$ in Equation \ref{eq:reduced_dissipation_inequality}.
Figure \ref{fig:reactionforces_gradient_hardening} displays the reaction force response of the three cases. All three scenarios first undergo a linear elastic loading phase, followed by hardening and finally softening behavior. The linear elastic behavior is obtained even though an AT2-type phase field model is applied, since the material degradation is based on a combination of accumulated plastic strain and the phase fields (compare Equation (\ref{eq:degradation_function})). As a consequence, the local phase field can only develop once at least a small amount of plastic strain develops, thus recovering a true linear elastic phase. As expected, the gradient hardening together with micro-hard conditions gives a stiffer response in the nonlinear regime. 
Figure \ref{fig:degradation_gradient_hardening} shows the material degradation at the last converged load step.
The well known diamond shape crack pattern (compare e.g. \cite{DeLorenzis2016}) is recovered without gradient hardening. The behavior under micro-free boundary conditions is similar, but shows smoothing of the degradation field. Micro-hard boundary conditions disallow the development of plastic strains on the boundary in a weak sense (for the chosen algorithmic formulation). Consequently, no damage can develop on the boundaries for a sufficiently fine mesh.
The micro-hard response is more sensitive to the mesh, which is shown in Figure \ref{fig:degradation_gradient_hardening} by the fact that an unsymmetric mesh gives a slightly unsymmetric response. The mesh for the micro-hard case has however already been refined compared to the other cases by using elements with an average size of $0.25\,\text{mm}$ on the left and right boundaries.

\subsection{Effect of micromorphic penalty parameter}
\label{subsection:example_micromorphic}
\begin{figure}
    \centering    \includegraphics{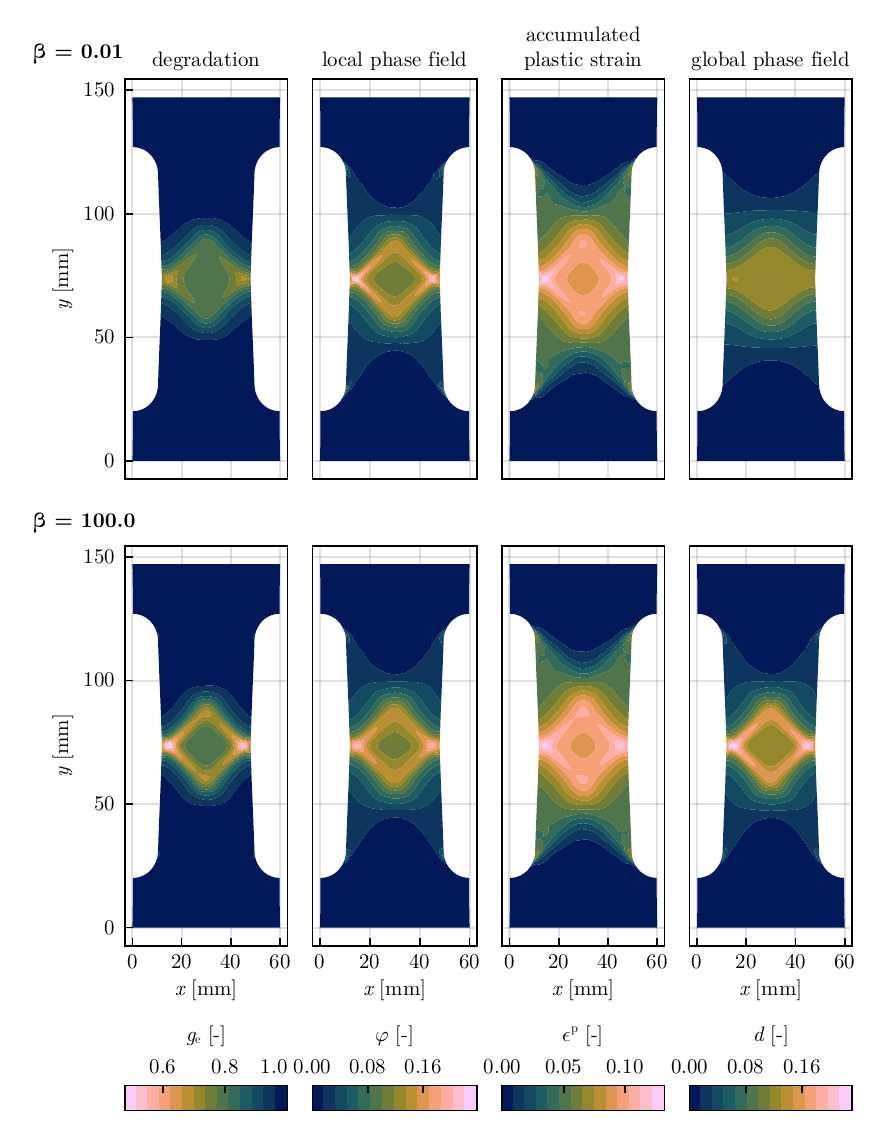}
    \caption{Effect of the micromorphic penalty parameter $\beta$: The top row shows results for a low penalty effect, the bottom row shows results for a sufficient penalty effect. If the penalty parameter is chosen too low a decoupling between the global and the local phase field occurs, resulting in a \textit{lack} of regularization of the local phase field $\varphi$. In this case the local phase field $\varphi$ is almost exclusively driven by the the accumulated plastic strain $\epsilon^\mathrm{p}$, which leads to a high level of localization in the local phase field. The global phase field $d$ instead experiences a lack of localization, as the coupling to the local phase field $\varphi$ decreases.}
    \label{fig:micromorphic_parameter}
\end{figure}
The micromorphic phase-field formulation adopts a penalty parameter $\alpha = \beta \, \mathcal{G}^\mathrm{d}_0/\ell_0$. \citeauthor{Bharali2023} \cite{Bharali2023} have shown for linear elastic problems that insufficient penalization leads to premature failure of the specimen (since the formulation turns into a local damage model in this case). For the chosen model problem, we have observed that the reaction force response is relatively robust against the choice of penalty parameter, but shows the same tendency of premature failure for very low choices of $\beta$. 
The lower sensitivity in our case can largely be attributed to the influence of plasticity in the proposed model. Figure \ref{fig:micromorphic_parameter} displays the effect of drastically reducing the micromorphic penalty parameter $\beta$. Micro-free boundary conditions have been applied in these simulations. 
The resulting effect of a lower $\beta$ is a decoupling of the local phase field $\varphi$ and the global phase field $d$ (as expected, compare Equation (\ref{eq:free_energy_micromorphic})). While both phase fields show a similar shape and magnitude for $\beta=100$, the global phase field $d$ has much lower values and lacks localization for $\beta=0.01$. In contrast, the local phase field over-localizes for an insufficient penalty parameter in this example. The reason for this lack of regularization is that the local phase field $\varphi$ develops due to the increase of accumulated plastic strain $\epsilon^\mathrm{p}$ in this case, it is mostly uninfluenced by the global phase field $d$. This means that the effect of crack regularization inherent to phase field fracture modeling, which occurs in the global equation system (compare Equation (\ref{eq:weak_pf})), does not permeate to the local problem described in Equation (\ref{eq:local_residual_pf}), which determines $\varphi$.
It should be noted that gradient-enhanced plasticity represents a regularization of the plastic strain field and that the response for an unsufficiently penalized micromorphic phase-field model without gradient-enhanced plasticity shows even more pathological localization.

\subsection{Damage irreversibility}
\label{subsection:example_cyclic}
\begin{figure}[htb]
    \centering
    \includegraphics{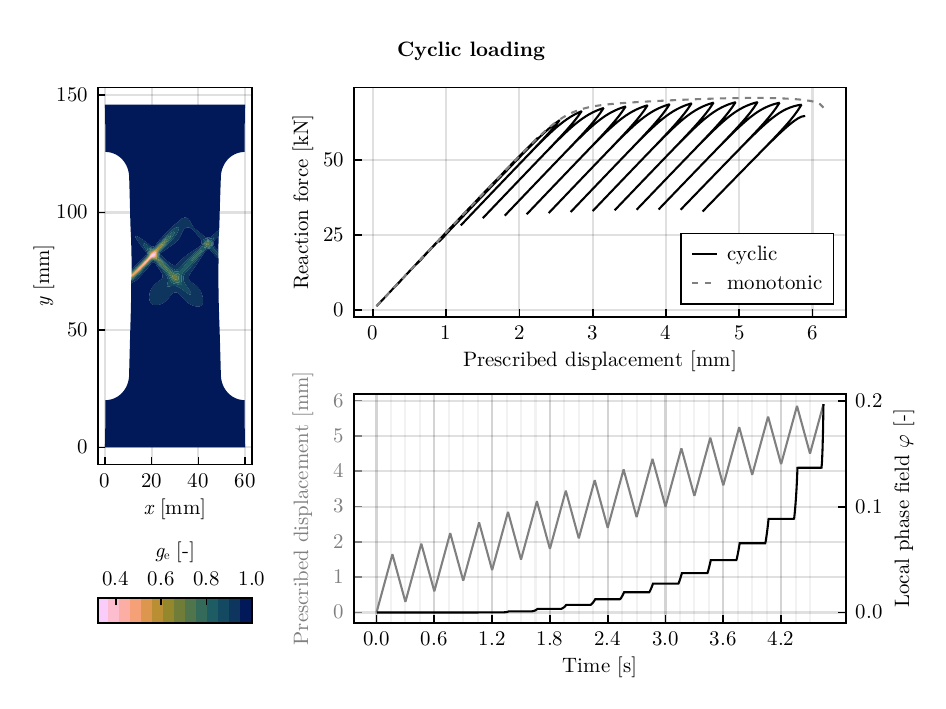}
    \caption{Damage irreversibility can be observed under ratcheting loading. The I-shaped specimen with inhomogeneities is loaded under racheting loading in displacement control. Elastic unloading and reloading with a visible viscous impact can be observed in the reaction force response (top right). The prescribed loading curve and the development of the local phase field $\varphi$ in a specific material point are shown on the bottom right. After the initial elastic phase the local phase field develops under advancing loading, but remains constant under unloading and reloading. The left plot shows the crack shape at the end of the simulation.}
    \label{fig:cyclic}
\end{figure}

In the presented model, damage irreversibility is ensured by introducing a history variable for the local phase field, compare Section \ref{section:irreversibility}. In order to demonstrate the behavior of this formulation, the I-shaped specimen is loaded in a ratcheting manner, whereby the loading rate of $10\,\mathrm{mm/s}$ is kept from the previous numerical example. The loading curve is shown in the bottom right panel of Figure \ref{fig:cyclic}.
In addition to the previous setup, two inhomogeneities are added to the I-shaped specimen at positions $\boldsymbol{X}_\mathrm{c1}=\left[45.0\,\mathrm{mm}, 83.125\,\mathrm{mm}\right]$ and $\boldsymbol{X}_{\mathrm{c}2}=\left[20.0\,\mathrm{mm}, 78.75\,\mathrm{mm}\right]$.
The inhomogeneities are introduced by smoothly reducing the yield limit in the vicinity of the points $\boldsymbol{X}_{\mathrm{c}i}$ by up to $95\,\%$, such that the effective yield limit is
\begin{equation}
    \label{eq:inhomogeneity}
    \tau_{\mathrm{y,red}}
    =
    \left(
        1 - 0.95\,\mathrm{b}
        \left(
            \frac{\left|\boldsymbol{X} - \boldsymbol{X}_{\mathrm{c}i} \right|}{r_\mathrm{red}}
        \right)
    \right)
    \tau_\mathrm{y} \,,
\end{equation}
where $\mathrm{b}$ is a bump function
\begin{equation}
    \label{eq:bump}
    \mathrm{b}\left(x\right)
    =
    1 - \frac{\exp\left(-1/x\right)}{\exp\left(-1/x\right)+\exp\left(-1/(1-x)\right)} \,.
\end{equation}
A radius $r_\mathrm{red}$ of $4\,\mathrm{mm}$ is chosen for the inhomogeneities and the mesh is refined with an average element size of $0.25\,\mathrm{mm}$ within this radius around the points $\boldsymbol{X}_{\mathrm{c}i}$. 
The left part of Figure \ref{fig:cyclic} shows the degradation $g_\mathrm{e}$ at the final step of the cyclic simulation. On the upper right of the figure, the resulting reaction forces from the cyclic loading and the equivalent monotonic loading are shown. On the lower right of the figure, the prescribed displacement together with the development of the local phase field $\varphi$ in a heavily degraded material point are shown. 
Within the first two cycles elastic loading and unloading are observed. During the subsequent cycles damage starts to develop in the vicinity of the inhomogeneities, which has an impact on the global reaction force curve. A noticeable viscous impact can be observed from the reaction force response around the load reversals and when comparing to the reaction force resulting from monotonic loading. The local phase field history in the chosen integration point is displayed in the bottom right graph of Figure \ref{fig:cyclic}. It begins to develop after 4 cycles. The phase field grows during progressive loading, but remains constant during unloading and reloading. The same behavior can be observed for the degradation $g_\mathrm{e}$ and the global phase field $d$.
The obtained crack pattern shows a clear influence of the circular inhomogeneities. The diamond-shaped crack pattern observed in the previous numerical examples is disrupted and instead a zig-zag shaped crack pattern results from the slip system orientation in conjunction with the inhomogeneities.

\subsection{3D simulation}
\label{subsection:example_3d}
\begin{figure}[htb]
    \centering
    \includegraphics[width=15cm]{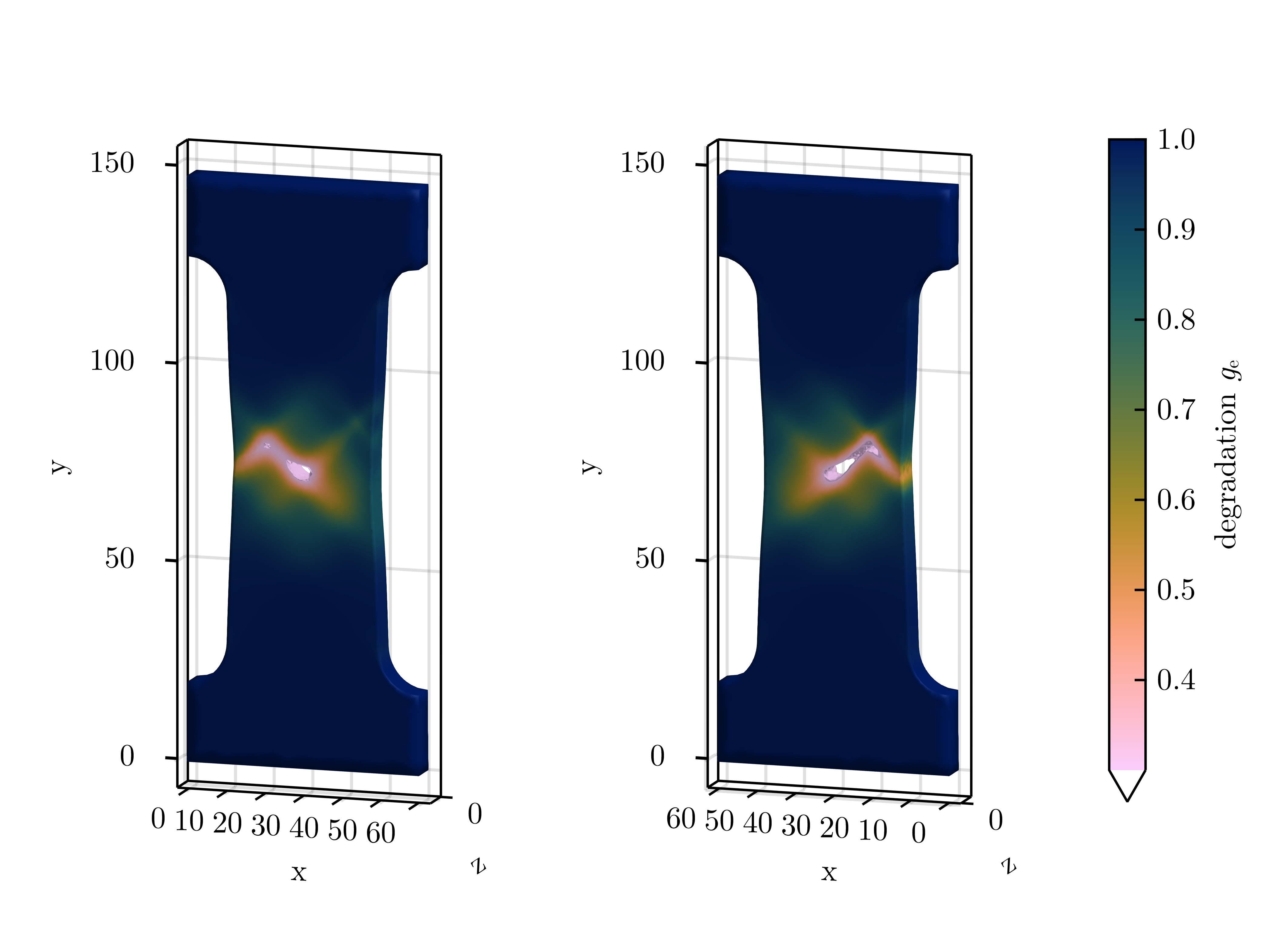}
    \caption{Degradation $g_\mathrm{e}$ in the three-dimensional I-shaped specimen with two spherical inclusions at different through-thickness coordinates. Degradation values below 0.3 are displayed as crack. Distinct crack patterns evolve on the front (left plot) and back (right plot) side of the specimen, highlighting the three-dimensional nature of the crack.
    All axes are measured in mm. }
    \label{fig:3d_degradation}
\end{figure}

A final simulation is conducted in order to display the capability of the model and numerical implementation to capture crack fronts in three dimensional space. The I-shaped specimen is extruded to a thickness of $8\,\mathrm{mm}$. Following Equations (\ref{eq:inhomogeneity}) and (\ref{eq:bump}), spherical inhomogeneities are introduced at positions 
$\boldsymbol{X}_\mathrm{c1}=\left[45.0\,\mathrm{mm}, 83.125\,\mathrm{mm}, 5.0\,\mathrm{mm}\right]$ and $\boldsymbol{X}_{\mathrm{c}2}=\left[20.0\,\mathrm{mm}, 78.75\,\mathrm{mm}, 3.0\mathrm{mm}\right]$. Thereby a through-thickness inhomogeneous response is expected. Micro-free boundary conditions are applied and the phase field length scale $\ell_0$ is set to $2\,\mathrm{mm}$. 
Figure \ref{fig:3d_degradation} displays the degradation $g_\mathrm{e}$ at the final step of the simulation from the front and the back side of the sample. For degradation values below $0.3$, the material is assumed to be fractured and is eliminated from the results. While the overall degradation response is comparable to that from the similar 2D cyclic simulation, differences in the crack patterns on the front and back side of the sample can be observed. A three-dimensional crack front develops, resulting from the two inhomogeneities placed at different through-thickness positions.
\section{Concluding remarks}
\label{section:concluding}
We have presented a thermodynamical framework that incorporates gradient-enhanced crystal plasticity in conjunction with a ductile phase-field fracture model in a large deformation setting.
While the framework allows different approaches to incorporate damage irreversibility, this work puts the major focus on testing a recently suggested variationally and thermodynamically consistent micromorphic formulation, that allows for a localized phase-field formulation and enforces the irreversibility constraint on an integration point level.
A prototype model is introduced for performing numerical experiments of the presented framework. It adopts an AT2-type of phase-field model which gives ductile fracture behavior by the choice of degradation function.   
The gradient-enhanced crystal plasticity model in conjunction with the micromorphic phase-field formulation and the ductile (plastic-strain based) degradation function leads to local equation systems that involve full coupling between local variables pertaining to the equilibrium equations as well as to the phase field equation. 
Therefore, a staggered solution scheme that extends the global staggered solution scheme to the local level is presented and adopted.
In the numerical examples, we demonstrate the effect that gradient-enhanced plasticity and the associated boundary conditions have on the fracture response.
It is shown how the micromorphic phase-field formulation relies on sufficient penalization of the difference between the global and the local phase fields. For unsufficient penalization, the model turns into a local damage model and suffers from the associated drawbacks.
Further, our numerical examples display that the presented model predicts unloading and reloading in a physically meaningful manner. The degradation of the material remains constant during unloading and advances once reloading goes past the previous level.
Finally, the cyclic and the 3D experiments show that the model is able to account for material inhomogeneities and the resulting (arbitrary) crack patterns.

The presented work opens up for future investigation of polycrystalline fracture, where the authors in particular are interested in investigating the effect of grain boundaries in polycrystalline transgranular fracture. The presented gradient-extended crystal plasticity model allows to study the interaction between slip transfer at grain boundaries and crack propagation, as well as the impact of size effects in crack formation on the micro-scale. In order to study crack propagation, solution schemes that allow to trace unstable crack growth, such as arc-length schemes (\cite{Bharali2022}, \cite{Borjesson2022}), and strategies to address convergence problems in the local equations systems are deemed important.

\FloatBarrier
\section*{Acknowledgements}
The work in this paper has been funded by the Swedish Research Council (Vetenskapsrådet) under the grant number 2018-04318.

\FloatBarrier
\newpage
\printbibliography

\newpage
\FloatBarrier
\begin{appendices}
\section{Material parameters}
\label{section:appendix_parameters}
\begin{table}[h!]
    \centering
{\linespread{1.1}\selectfont
    \begin{tabular}{||l l r l||}
        \hline
        Parameter &  & Value & Unit \\ [0.5ex]
        \hline \hline
        Bulk modulus & $\kappa$ & 71660 & MPa \\
        \hline
        Shear modulus & $\mu$ & 27260  & MPa \\ [0.5ex]
        \hline \hline
        Yield stress & $\tau^\mathrm{y}$ & 345 & MPa \\ [0.5ex]
        \hline
        Isotropic hardening modulus & $H_\mathrm{iso}$ & 250 & MPa \\
        \hline
        Gradient hardening modulus & $H^\mathrm{g}$ & 1000 & MPa \\
        \hline
        Gradient hardening length scale & $l^\mathrm{g}$ & 4.0 & mm \\
        \hline \hline
        Visco-plastic relaxation time & $t^\ast$ & 1 & s \\
        \hline
        Visco-plastic drag stress & $\sigma^\mathrm{d}$ & 500 & MPa \\
        \hline
        Visco-plastic exponent & $m$ & 8 & - \\ [0.5ex]
        \hline \hline
        Effective fracture energy & $\mathcal{G}^\mathrm{d}_0 / \ell_0$ & 300 & $\mathrm{N}/\mathrm{mm}^2$ \\
        \hline
        Phase field length scale & $\ell_0$ & 0.5 & mm \\
        \hline
        Micromorphic interaction parameter & $\alpha$ & 60000 & $\mathrm{N/mm^2}$ \\
        \hline
        Critical plastic strain & $\epsilon^\mathrm{p}_\mathrm{crit}$ & 0.1 & - \\
        \hline
        Degradation exponent & $n$ & 2 & - \\
        \hline
    \end{tabular}
    }
    \caption{Base material parameters employed for the numerical experiments.}
    \label{tab:material_parameters}
\end{table}

\begin{table}[h!]
    \centering
{\linespread{1.2}\selectfont
    \begin{tabular}{||c c c|c c c|c c c|c c c||}
        \hline
       $\alpha$  & $\bar{\boldsymbol{s}}_\alpha$ & $\bar{\boldsymbol{m}}_\alpha$ &
       $\alpha$  & $\bar{\boldsymbol{s}}_\alpha$ & $\bar{\boldsymbol{m}}_\alpha$ & 
       $\alpha$  & $\bar{\boldsymbol{s}}_\alpha$ & $\bar{\boldsymbol{m}}_\alpha$ &
       $\alpha$  & $\bar{\boldsymbol{s}}_\alpha$ & $\bar{\boldsymbol{m}}_\alpha$ \\ [0.5ex]
       \hline \hline
       1 & $\left[\bar{1} \, 1 \, 0 \right]$ & $\left[1 \, 1 \, 1 \right]$ &
       4 & $\left[\bar{1} \, \bar{1} \, 0 \right]$ & $\left[1 \, \bar{1} \, \bar{1} \right]$ &
       7 & $\left[1 \, 1 \, 0 \right]$ & $\left[\bar{1} \, 1 \, \bar{1} \right]$ &
       10 & $\left[1 \, \bar{1} \, 0 \right]$ & $\left[\bar{1} \, \bar{1} \, 1 \right]$ \\
       \hline
       2 & $\left[1 \, 0 \, \bar{1} \right]$ & $\left[1 \, 1 \, 1 \right]$ &
       5 & $\left[1 \, 0 \, 1 \right]$ & $\left[1 \, \bar{1} \, \bar{1} \right]$ &
       8 & $\left[\bar{1} \, 0 \, 1 \right]$ & $\left[\bar{1} \, 1 \, \bar{1} \right]$ &
       11 & $\left[\bar{1} \, 0 \, \bar{1} \right]$ & $\left[\bar{1} \, \bar{1} \, 1 \right]$ \\
       \hline
       3 & $\left[0 \, \bar{1} \, 1 \right]$ & $\left[1 \, 1 \, 1 \right]$ &
       6 & $\left[0 \, 1 \, \bar{1} \right]$ & $\left[1 \, \bar{1} \, \bar{1} \right]$ &
       9 & $\left[0 \, \bar{1} \, \bar{1} \right]$ & $\left[\bar{1} \, 1 \, \bar{1} \right]$ &
       12 & $\left[0 \, 1 \, 1 \right]$ & $\left[\bar{1} \, \bar{1} \, 1 \right]$ \\
       \hline
    \end{tabular}
}
    \caption{FCC slip systems. The unit cell is aligned with the coordinate axes in the numerical examples. In the numerical examples slip systems $2,\,5,\,8,\,11$ are omitted, because they are perpendicular to the loading loading direction and have a negligible impact on the plastification process. In the cyclic and the 3D example, only the 4 slip systems with the largest impact are used, that is slip systems $1,\,4,\,7,\,10$.}
    \label{tab:slipsystems}
\end{table}
\FloatBarrier
\section{Microforce derivation}
\label{section:microforce}
An alternative derivation of the gradient extended plasticity and the phase-field following the procedure in \cite{Borden2018} is summarized here in point-form. 

\begin{itemize}
    \item Microforce balances
\begin{eqnarray}
   \bi \nabla_0 \cdot \bi \xi_\alpha+\pi_\alpha+l_\alpha&=&0 \\
    \bi \nabla_0 \cdot \bi \xi_d+\pi_d+l_d&=&0
\end{eqnarray}
with $l_\alpha=l_d=0$.
\item Dissipation inequality
\begin{equation}
    \bi P : \dot{\bi F}+\bi \xi_\alpha \cdot \nabla_0 \dot{k}_\alpha-\pi_\alpha \, \dot{k}_\alpha+
    \bi \xi_d \cdot \nabla_0 \dot{d}-\pi_d \, \dot{d}-\dot{\Psi} \geq 0
\end{equation}
\item Free energy
\begin{equation}
    \Psi=\Psi\left(\bi C_{\rm e}, \, q, \, k_\alpha, \, \boldsymbol{\nabla}_0 \, k_\alpha, \, d, \, \boldsymbol{\nabla}_0 d \, \right)
\end{equation}
\item gives
\begin{equation}
    \bi S_{\rm e}=2 \, \frac{\partial \Psi}{\partial \bi C_{\rm e}}
    \label{eq:elastic_2nd_PK_stress_from_microforce}
\end{equation}
and the reduced dissipation inequality
\begin{equation}
\mathcal{D} =
      {\boldsymbol{M}}_\mathrm{e}  
                :
                {\boldsymbol{L}}_\mathrm{p}  
                +Q \, \dot{q}
    -\left( 
        \frac{\partial\Psi}{\partial k_\alpha} +\pi_\alpha
 \right)\, \dot{k}_\alpha
    -
 \left( 
        \frac{\partial\Psi}{\partial \boldsymbol{\nabla}_0 k_\alpha} 
        -\bi \xi_\alpha
        \right) 
        \cdot \boldsymbol{\nabla}_0 \dot{k}_\alpha
    -\left( 
        \frac{\partial\Psi}{\partial d} +\pi_d
 \right)\, \dot{d}
    -
 \left( 
        \frac{\partial\Psi}{\partial \boldsymbol{\nabla}_0 d} 
        -\bi \xi_d
        \right) 
        \cdot \boldsymbol{\nabla}_0 \dot{d}
    \geq
    0
\end{equation}
\item By assuming  $\bi \xi_\alpha$, $\pi_d$ and $\bi \xi_d$ energetic, i.e.
\begin{eqnarray}
%\pi_\alpha&=&-  \frac{\partial\Psi}{\partial k_\alpha}, \quad
\bi \xi_\alpha= \frac{\partial\Psi}{\partial \boldsymbol{\nabla}_0 k_\alpha}, \quad
\pi_d=-  \frac{\partial\Psi}{\partial d}, \quad
\bi \xi_d= \frac{\partial\Psi}{\partial \boldsymbol{\nabla}_0 d} 
\end{eqnarray}
gives
\begin{equation}
      {\boldsymbol{M}}_\mathrm{e}  
                :
                {\boldsymbol{L}}_\mathrm{p}  
                +Q \, \dot{q}
    +\kappa_\alpha
    %-\left( 
     %   \frac{\partial\Psi}{\partial k_\alpha} +\pi_\alpha
 %\right)
 \, \dot{k}_\alpha
    \geq
    0
\end{equation}
where
\begin{equation}
    \kappa_\alpha=-
        \frac{\partial\Psi}{\partial k_\alpha} +\bi \nabla_0 \cdot \frac{\partial\Psi}{\partial \boldsymbol{\nabla}_0 k_\alpha}
\end{equation}
by using the microforce balance. This is identical to when $Y_\alpha=0$ in (\ref{eq:kappa_def}).
\item Insert $\pi_d$ and $\bi \xi_d$ into microforce balance
\begin{eqnarray}
    \bi \nabla_0 \cdot\frac{\partial\Psi}{\partial \boldsymbol{\nabla}_0 d} -\frac{\partial\Psi}{\partial d}=0
\end{eqnarray}
which is identical to $Y=0$ in (\ref{eq:Y_def}).

\end{itemize}

\end{appendices}

\end{document}